\documentclass[12pt,a4paper]{article}
\usepackage{amsmath,amssymb,amsthm,fullpage,color}
\usepackage{hyperref}
\hypersetup{colorlinks=true,linkcolor=blue,citecolor=red}
\newtheorem{theorem}{Theorem}[section]
\newtheorem{thm}[theorem]{Theorem}
\newtheorem{prop}[theorem]{Proposition}
\newtheorem{lem}[theorem]{Lemma}
\newtheorem{cor}[theorem]{Corollary}

\numberwithin{equation}{section}

\newlength{\bibitemsep}
\newlength{\bibparskip}
\let\oldthebibliography\thebibliography
\renewcommand\thebibliography[1]{%
  \oldthebibliography{#1}%
  \setlength{\parskip}{\bibitemsep}%
  \setlength{\itemsep}{\bibparskip}%
}
\def\d{\mathrm{d}}
\def\ind{\mathbf{1}}
\def\Nb{\mathbb{N}}
\def\Pb{\mathbb{P}}
\def\Eb{\mathbb{E}}
\def\Rb{\mathbb{R}}
\def\Zb{\mathbb{Z}}
\def\Ic{\mathcal{I}}
\def\Hc{\mathcal{H}}
\def\Pc{\mathcal{P}}

\def\Cc{\mathcal{C}}
\def\Tc{\mathcal{T}}
\def\Sc{\mathcal{S}}
\def\Pr{\mathbf{P}}
\def\Dk{\mathfrak{D}}
\def\BK{\mathcal{BK}}

\begin{document}

\title{Biased random walk on supercritical percolation: Anomalous fluctuations in the ballistic regime}
\author{A.~M.~Bowditch\footnote{School of Mathematics and Statistics, University College Dublin, Dublin 4, Ireland. Email: adam.bowditch@ucd.ie.} and D.~A.~Croydon\footnote{Research Institute for Mathematical Sciences, Kyoto University, Kyoto 606-8502, Japan. Email: croydon@kurims.kyoto-u.ac.jp.}}
%60K37   	Processes in random environments
%60G50   	Sums of independent random variables; random walks
%60K35   	Interacting random processes; statistical mechanics type models; percolation theory

\maketitle

\begin{abstract}
We study biased random walk on the infinite connected component of supercritical percolation on the integer lattice $\mathbb{Z}^d$ for $d\geq 2$. For this model, Fribergh and Hammond showed the existence of an exponent $\gamma$ such that: for $\gamma<1$, the random walk is sub-ballistic (i.e.\ has zero velocity asymptotically), with polynomial escape rate described by $\gamma$; whereas for $\gamma>1$, the random walk is ballistic, with non-zero speed in the direction of the bias. They moreover established, under the usual diffusive scaling about the mean distance travelled by the random walk in the direction of the bias, a central limit theorem when $\gamma>2$. In this article, we explain how Fribergh and Hammond's percolation estimates further allow it to be established that for $\gamma\in(1,2)$ the fluctuations about the mean are of an anomalous polynomial order, with exponent given by $\gamma^{-1}$.
\end{abstract}

\section{Introduction}

Understanding the interplay between the geometry of a space and the stochastic processes that live upon it has been a major focus of probability over the last four decades. Of particular interest is the case when the space in question is random, with a central example being a percolation cluster on the integer lattice $\mathbb{Z}^d$. Specifically, it is now well-known that for the simple random walk on the unique infinite cluster of supercritical percolation on $\mathbb{Z}^d$, one sees under the usual diffusive scaling Brownian motion as a scaling limit (see \cite{Bar,BB,DFGW,MP,SS} for some key works in this direction). Whilst this is the same qualitative behaviour as for simple random walk on the integer lattice itself, there is additional complexity in that the inhomogeneity of the environment incorporates traps, that is, areas such that if the random walk enters them, then it will take an anomalously long time to escape, and these lead to the simple random walk on the percolation cluster having a smaller diffusion constant than for the random walk on the whole lattice. In the case of critical percolation, it has been shown that such an effect is strong enough to lead to sub-diffusive behaviour of the associated simple random walk \cite{K}.

The situation of trapping is exacerbated by adding a bias to the random walk, so that it has a particular directional preference. Indeed, it has been observed for various models of random graph that trapping can become stronger when the bias is increased, with weak biases leading to ballistic (positive speed) behaviour asymptotically (as on $\mathbb{Z}^d$), but strong biases leading to sub-ballisticity. See \cite{BGP,F1,FH,S} for relevant mathematical work on supercritical percolation clusters, and \cite{Dhar} for an earlier physics discussion of the issue. We highlight that it was in \cite{FH} that the sharpness of this phase transition was observed. Moreover, as we will explain in more detail below, \cite{FH} also contained a description of the polynomial escape exponent in the sub-ballistic regime, and established, in a regime of suitably weak biases, a functional central limit theorem, demonstrating Gaussian fluctuations around the ballistic mean behaviour. For context, we note that these developments build on the substantial literature concerning the simpler setting of biased random walk on Galton-Watson trees. See \cite{A,AR,BFGH,BFS,BH,B1,B2,B4,CFK,H,LPP} for example works in the latter area, and \cite{BF} for a relatively recent and comprehensive survey of biased random walks on random graphs.

The focus of this article is on the regime where the biased random walk on the supercritical percolation cluster is ballistic, but expected to show non-Gaussian fluctuations around this, or more precisely, $\alpha$-stable fluctuations for some $\alpha\in (1,2)$. (See \cite[Conjecture 4.2]{BF} or the comment below \cite[Theorem 1.4]{FH}.) Understanding the fluctuations in this case is slightly more delicate than the results discussed in the previous paragraph since it concerns not only the leading behaviour, but also the second order asymptotics, of the process in question. To date, corresponding  $\alpha$-stable scaling limits have been observed for one-dimensional random walk in random environment \cite{KKS}, a model of one-dimensional model of randomly trapped random walk \cite{B3,B4}, and also for subcritical Galton-Watson trees conditioned to survive (also \cite{B4}). Whilst we do not prove such exact results for supercritical percolation, we will nonetheless show that biased random walk on supercritical percolation in the regime of interest exhibits fluctuations on the same polynomial scale.

We now proceed to introduce the model of interest in detail, starting with the random environment. Denote by $E(\Zb^d)$ the set of nearest-neighbour edges of the lattice $\Zb^d$ for $d\geq 2$. Fix $p\in(0,1)$, and let
\begin{equation*}
P_p=(p\delta_1+(1-p)\delta_0)^{\otimes E(\Zb^d)}
\end{equation*}
be the probability measure on $\{0,1\}^{E(\Zb^d)}$ corresponding to Bernoulli bond percolation on $\Zb^d$. (Here, $\delta_i$ represents the probability measure on $\{0,1\}$ placing all of its mass at $i$.) We say that an edge $e\in E(\Zb^d)$ is {\it open} in the configuration $\omega\in\{0,1\}^{E(\Zb^d)}$ if $\omega(e)=1$, and {\it closed} if $\omega(e)=0$. The set of open edges induces a subgraph of $(\Zb^d,E(\Zb^d))$, which we also denote $\omega$. It is classical that there exists a critical percolation probability $p_c(d)\in(0,1)$ such that there exists a unique infinite open cluster $P_p$-a.s.\ if $p>p_c(d)$, whereas no infinite cluster exists for $p<p_c(d)$. We henceforth suppose $p>p_c(d)$, and define
\[\Pr_p(\cdot)=P_p(\cdot \:\vert\: \Ic),\]
where $\Ic$ is the event that the origin is contained in the unique infinite open cluster.

We next define the biased random walk associated with a particular realisation of $\omega\in\Omega$. Firstly, we suppose the bias is given by $\ell=\lambda \vec{\ell}\in\mathbb{R}^d\backslash\{0\}$, where $\lambda>0$ represents the strength of the bias, and $\vec{\ell}$ is a vector in the unit sphere $\mathbb{S}^{d-1}$ that determines its direction. We set the conductance between nearest neighbours $x$ and $y$ in $\Zb^d$ to be the quantity
\[c^\omega(x,y)=e^{(x+y)\cdot \ell} \omega(\{x,y\}).\]
We then consider the discrete-time Markov chain $(X_n)_{n\geq 0}$ on $\Zb^d$ with law $P_x^\omega$, which satisfies $P_x^\omega(X_0=x)=1$ and whose transition probabilities $p^\omega(x,y)$ are given by
\[p^\omega(x,y) = \begin{cases} 1, & \text{if } x=y \text{ and } \omega(\{x,z\})=0 \text{ for all } z\text{ such that }|z-x|=1, \\ \frac{c^\omega(x,y)}{\sum_{z:\:|z-x|=1}c^\omega(x,z)}, & \text{otherwise.} \end{cases}\]
We define the annealed law of the random walk on the supercritical percolation cluster to be the semi-direct product
\[\Pb_p(\cdot)=\int_{\Omega}  P_0^\omega(\cdot) \Pr_p(\d\omega).\]

We are now in a position to describe the main results of \cite{FH}. To do this, let us start by introducing the so-called backtrack function: for $x\in\Zb^d$,
\begin{equation}\label{bt}
\BK(x)=\begin{cases} 0, & \text{if } x\notin \mathcal{C}_\infty, \\ \min_{(p_x(i))_{i\geq 0}\in\Pc_x}\max_{i\geq 0}(x-p_x(i))\cdot \vec{\ell}, & \text{otherwise,} \end{cases}
\end{equation}
where $\Pc_x$ is the set of infinite, open, self-avoiding paths started from $x$, and $\mathcal{C}_\infty$ is the unique infinite cluster of our supercritical percolation model. This represents the smallest distance against the bias that a random walker will have to travel from $x$ in order to escape from the trap that contains it, where we think of traps in the environment as being the connected components of vertices where $\BK$ is strictly positive. By \cite[Proposition 1.1]{FH}, there exists a constant $\zeta=\zeta(p,\vec{\ell},d)\in(0,\infty)$ such that
\[\lim_{n\to\infty}n^{-1}\log\Pr_p(\BK(0)>n)=-\zeta;\]
this exponent governs the sizes of traps that are seen in the environment. Since it takes time of order $e^{2\lambda h}$ to exit a trap of depth $h$, the exponent $\zeta$ turns out to also be strongly linked to the distribution of trapping times experienced by the biased random walk, though for the latter purpose, it is more natural to consider the parameter:
\begin{equation}\label{gammadef}
\gamma=\frac{\zeta}{2\lambda}.
\end{equation}
(Cf.\ discussion of \cite[Section 1.4]{FH}.) Indeed, it is established in \cite[Theorem 1.2]{FH} that for $\gamma>1$ (the weak bias regime),
\begin{equation}\label{ballistic}
\lim_{n\to\infty}\frac{X_n}{n}=\vec{v}, \qquad  \Pb_p\text{-a.s.},
\end{equation}
where $\vec{v}\in\mathbb{R}^d$ is a deterministic vector satisfying $\vec{v}\cdot \ell>0$, whereas for $\gamma<1$ (the strong bias regime), the above limit holds with $\vec{v}=0$. (A phase transition was shown earlier in \cite{BGP,S}, but the change-point was not determined.) A rough intuitive explanation for this phenomenon is that to reach distance $n$ from the origin, the random walker has to pass through order $n$ approximately independent and identically distributed (i.i.d.) traps, and the time it takes to escape each of these has a distribution that decays polynomially with exponent $\gamma$. When $\gamma>1$, this distribution has a finite mean, and we thus see ballistic, law of large numbers-type behaviour; when $\gamma< 1$, the trapping time distribution has an infinite mean, which leads to sub-ballistic behaviour. In the latter case, in \cite[Theorem 1.5]{FH} it is clarified that, $\mathbb{P}_p$-a.s.,
\begin{equation}\label{gamma1}
\lim_{s\to\infty}\frac{\log \Delta_s}{\log s}=\frac{1}{\gamma},
\end{equation}
where
\[\Delta_s=\inf\{m\in\Nb:X_m\cdot\vec{\ell}\geq s\},\]
is the first time that the random walk exceeds distance $s$ in direction $\vec{\ell}$, and also
\begin{equation}\label{gamma2}
\lim_{n\to\infty}\frac{\log X_n\cdot\vec{\ell}}{\log n}=\gamma,
\end{equation}
as might be expected from the previous discussion. (Note that, in contrast to the discrete parameter $n$ of \eqref{gamma2}, the parameter $s$ of \eqref{gamma1} is continuous.) Moreover, when $\gamma>2$ -- the case when the trapping time distribution has finite variance, it is shown as \cite[Theorem 1.4]{FH} that the model exhibits central-limit theorem-type behaviour, in that
\[\left(\frac{X_{tn}-tn\vec{v}}{\sqrt{n}}\right)_{t\geq 0}\]
converges under the annealed law $\mathbb{P}_p$ to a non-degenerate Brownian motion.

Our aim is to fill in the part of the story concerning the scale of the fluctuations around \eqref{ballistic} in the non-Gaussian case. In particular, we prove the following, where we define $\log_+(x)=\max\{0,\log(x)\}$.

\begin{thm}\label{mainthm} Let $d\geq 2$, $p\in(p_c(d),1)$, $\vec{\ell}\in\mathbb{S}^{d-1}$, $\gamma$ be defined as at \eqref{gammadef} and $v=\vec{v}\cdot\vec{\ell}$. If $\gamma\in(1,2)$, then
\begin{equation}\label{deltalim}
\frac{\log \left| \Delta_s-sv^{-1}\right|}{\log s}\buildrel{\mathbb{P}_p}\over\rightarrow\frac{1}{\gamma},
\end{equation}
and
\begin{equation}\label{xlim}
\frac{\log \left| X_n\cdot\vec{\ell}-nv\right|}{\log n}\buildrel{\mathbb{P}_p}\over\rightarrow\frac{1}{\gamma}.
\end{equation}
Moreover, it $\mathbb{P}_p$-a.s.\ holds that
\begin{equation}\label{limsup}
\limsup_{s\rightarrow\infty}\frac{\log \left| \Delta_s-sv^{-1}\right|}{\log s}=\limsup_{n\rightarrow\infty}\frac{\log \left| X_n\cdot\vec{\ell}-nv\right|}{\log n}=\frac{1}{\gamma}
\end{equation}
and
\begin{equation}\label{liminf}
\liminf_{s\rightarrow\infty}\frac{\log_+ \left| \Delta_s-sv^{-1}\right|}{\log s}=\liminf_{n\rightarrow\infty}\frac{\log_+ \left| X_n\cdot\vec{\ell}-nv\right|}{\log n}=0.
\end{equation}
\end{thm}

To check the above statements, we will use the regeneration structure of \cite{FH,S}, which will allow a useful comparison with a sequence of i.i.d.\ random variables. For these random variables, an upper distributional tail bound that will be sufficient for our purposes was obtained in \cite{FH}. Moreover, we adapt the arguments of \cite{FH} to give a corresponding lower bound. Although the leading order behaviour of the upper and lower bounds is polynomial with the same exponent, the possibility is left open that they differ by a sub-polynomial expression, and so we can not immediately apply classical results concerning i.i.d.\ sums. To overcome this, we develop machinery that allows the treatment of what might be considered `near-stable' random variables. As a further comment on the content of Theorem \ref{mainthm}, note that the relevant exponent only appears as a limit in $\mathbb{P}_p$-probability, rather than $\mathbb{P}$-a.s., as it did in \eqref{gamma1} and \eqref{gamma2}. This results from the oscillations around the mean that occur $\mathbb{P}$-a.s., which explain the differences between the limsups and liminfs of \eqref{limsup} and \eqref{liminf}.

Finally, it is natural to conjecture that in the settings where \eqref{gamma1}, \eqref{gamma2} and Theorem \ref{mainthm} were obtained, it would be possible to check distributional convergence results, with $\gamma$-stable random variables appearing in the limit, cf.\ \cite{B4,KKS}. (Precisely, due to lattice effects, the limits might only exist subsequentially; see \cite[Conjecture 4.2]{BF} and \cite{GK} for some more detailed discussion in this direction.) Other models where a similar trapping regime is expected to be found (and hence to which the results of this article might also apply) include biased random walk on supercritical Galton-Watson trees \cite{B1}, in random conductances \cite{F,FK}, on the interlacement set \cite{FP}, and on the trace of another biased random walk \cite{CH2,CH1}. One might also expect to see similar behaviour for random walk in transient random environments for which the tail of a suitable regeneration time distribution is of a suitable form, cf. \cite{S1}. We provide a more specific comment concerning the application of the results of this article to other models at the start of Section \ref{spsec}.

The remainder of the article is organised as follows. In Section \ref{nssec}, we establish fluctuation results for general near-stable i.i.d.\ random variables. In Section \ref{spsec}, we apply these results to the model of biased random walk on supercritical percolation. We remark that, unless otherwise noted, constants of the form $c$ and $C$ are deterministic, take values in $(0,\infty)$, and may change value from line to line.

\section{Fluctuations of near-stable random variables}\label{nssec}

Motivated by the corresponding problem for biased random walk on supercritical percolation clusters, in this section, we derive results concerning the fluctuations around a centring process for the partial sums of i.i.d.\ random variables. Throughout, we write $(\xi_i)_{i\geq1}$ for such a sequence, $\xi$ for a generic element of the sequence and $S_n=\sum_{i=1}^n\xi_i$. Our particular interest is in the case when $\xi$ has a distribution that is suitably close to an $\alpha$-stable distribution for some $\alpha\in(1,2)$. We present our `in probability' and `almost-sure' results separately. In the preparatory conclusions of Proposition \ref{p:Upper} and \ref{p:Lower}, we give a finer level of detail than is needed for the application of interest in this paper.

\subsection{Asymptotics in probability}

In order to state our main conclusions, we start by recalling that a (positive measurable) function $L$ is called slowly varying at $\infty$ if for any $s>0$ we have that
\[\lim_{t\rightarrow \infty}\frac{L(ts)}{L(t)}=1.\]
Moreover, a (positive measurable) function $f$ is called regularly varying with index $\alpha$ if it can be written as $f(t)=t^\alpha L(t)$ for a slowly varying function $L$. By \cite[Theorem 1.5.12]{BGT}, if a function $f$ is regularly varying with index $\alpha>0$, then there exists a (positive measurable) function $g$ called the asymptotic inverse of $f$ that is regularly varying with index $1/\alpha$ and satisfies
\[f(g(t))\sim g(f(t)) \sim t;\]
the function $g$ is asymptotically unique and can be chosen to be non-decreasing.

The first proposition we give represents an upper bound on the fluctuations of $S_n$ about its mean under the condition that we have a regularly varying upper bound on the tail of the distribution of $|\xi|$, with index $\alpha\in (1,2)$. We postpone the proofs of the all the results we state until later in the section.

\begin{prop}\label{p:Upper}
Let $\alpha\in(1,2)$ and $L$ be slowly varying at $\infty$ such that it is bounded away from $0$ and $\infty$ on every compact interval of $[0,\infty)$. Suppose that
\begin{equation}\label{e:upp}
\mathbf{P}\left(|\xi|\geq t\right)\leq t^{-\alpha}L(t),\qquad\forall t\geq1,
\end{equation}
and let $g$ be any non-decreasing asymptotic inverse of $t^\alpha/L(t)$. It is then the case that: for any $\rho<\alpha$, we can choose $C_\rho$ such that
\[\mathbf{P}\left(\left|\sum_{i=1}^n\frac{\xi_i-\mathbf{E}[\xi_i]}{a_n}\right|>\lambda\right)\leq C_\rho\lambda^{-\rho},\]
uniformly over $n,\lambda\geq 1$, where $a_n=g(n)$. In particular, the sequence
\[\left(\left|\frac{S_n-n\mathbf{E}[\xi]}{a_n}\right|\right)_{n\geq 1}\]
is tight.
\end{prop}

Our second proposition gives the complementary result, in that it contains a lower bound on the fluctuations of $S_n$ about any centring sequence under the condition that we have a regularly varying lower bound on the tail of the distribution of $|\xi|$, with index $\alpha>0$. (Note that for this result, we do not require that $\alpha\in(1,2)$.)

\begin{prop}\label{p:Lower}
Let $\alpha>0$ and $L$ be slowly varying at $\infty$ such that it is bounded away from $0$ and $\infty$ on every compact interval of $[0,\infty)$. Suppose that
\begin{equation}\label{e:low}
\mathbf{P}\left(|\xi|\geq t\right)\geq t^{-\alpha}L(t),\qquad\forall t\geq1,
\end{equation}
and let $g$ be any non-decreasing asymptotic inverse of $t^\alpha/L(t)$. It is then the case that: for any deterministic sequence $(c_n)_{n\geq 1}$ and any $\rho<\alpha/4$, there exists a constant $C_\rho$ such that
\begin{equation*}
\mathbf{P}\left(\left|\frac{S_n-c_n}{a_n}\right|<\lambda\right)\leq C_\delta\lambda^{\rho}
\end{equation*}
uniformly over $n\geq 1$ and $\lambda\geq 1/a_n$, where $a_n=g(n)$. In particular, the sequence
\begin{equation*}
\left(\left|\frac{S_n-c_n}{a_n}\right|^{-1}\right)_{n\geq 1}
\end{equation*}
is tight.
\end{prop}

Putting the previous two propositions together, we arrive at the following result. In particular, this states precisely that, even if we only have tail bounds for the distribution of $\xi$ that are tight up to small polynomial errors, it is possible to deduce the same fluctuation exponent as for random variables in the domain of attraction of a stable distribution of the same index.

\begin{cor}\label{c:convProb}
Let $\gamma\in(1,2)$. Suppose that, for every $\varepsilon\in(0,1)$, there exists a constant $C\in(0,\infty)$ such that
\begin{equation}\label{upperandlower}
C^{-1}t^{-(1+\varepsilon)\gamma}\leq \mathbf{P}\left(|\xi|\geq t\right)\leq Ct^{-(1-\varepsilon)\gamma},\qquad\forall t\geq1.
\end{equation}
It is then the case that
\[\frac{\log \left| S_n-n\mathbf{E}\xi\right|}{\log n}\buildrel{\mathbf{P}}\over\rightarrow\frac{1}{\gamma}.\]
\end{cor}

For use in the proofs, we next state a technical bound concerning regularly varying functions that constitutes one of the statements of Potter's Theorem \cite[Theorem 1.5.6]{BGT}.

\begin{thm}[Potter's theorem]\label{t:Pot}
Suppose $L$ is slowly varying. If $L$ is bounded away from $0$ and $\infty$ on every compact interval of $[0,\infty)$, then for every $\delta>0$, there exists $A_\delta$ such that: for all $x,y>0$,
\[\frac{L(x)}{L(y)}\leq A_\delta\max\left\{\left(\frac{x}{y}\right)^\delta,\left(\frac{x}{y}\right)^{-\delta}\right\}.\]
\end{thm}

We further have the following basic observation, which permits us to consider non-decreasing regularly varying functions.

\begin{lem}\label{l:nond} Suppose $f$ is regularly varying with index $\alpha>0$. It is then possible to find an asymptotically equivalent regularly varying function of the same index that is non-decreasing.
\end{lem}
\begin{proof} Let $g$ be an asymptotic inverse of $f$, as given by \cite[Theorem 1.5.12]{BGT}. This is regularly varying with index $1/\alpha>0$, and so it admits a non-decreasing asymptotic inverse $f^*$, which is regularly varying of index $\alpha$. Since $f$ is also an asymptotic inverse of $g$, by the asymptotic uniqueness of asymptotic inverses, it must be the case that $f$ and $f^*$ are asymptotically equivalent.
\end{proof}

We now adapt standard Laplace transform techniques for studying the distributional tails of sums of independent random variables (e.g.\ \cite[Section III.4]{P}) to prove Proposition \ref{p:Upper}. The estimates for the various integrals that appear in the proof are obtained using classical ideas from the study of stable laws (and regularly varying functions), of the kind developed in \cite{BGT, F2, IL, P}, for example.
\begin{proof}[Proof of Proposition \ref{p:Upper}]
First, suppose that $\xi$ is non-negative.
By Lemma \ref{l:nond}, we may assume that $t^{\alpha}/L(t)$ is non-decreasing. Applying Markov's inequality and the fact that $(\xi_i)_{i\geq 1}$ are i.i.d., we have that
\begin{align*}
\mathbf{P}\left(\sum_{i=1}^n(\xi_i-\mathbf{E}[\xi_i])<-\lambda a_n\right)
& \leq \inf_{\theta>0}\mathbf{P}\left(\exp\left(-\theta\sum_{i=1}^n(\xi_i-\mathbf{E}[\xi_i])\right)>\exp\left(\theta\lambda a_n\right)\right)\\
& \leq \inf_{\theta>0} \mathbf{E}\left[\exp\left(-\theta \xi\right)\right]^n\exp\left(n\theta\mathbf{E}[\xi]-\theta\lambda a_n\right).
\end{align*}
We moreover have that
\begin{align*}
\mathbf{E}\left[\exp\left(-\theta \xi\right)\right]
& = \int_0^1\mathbf{P}\left(\exp\left(-\theta \xi\right)>t\right)\d t \\
& = \theta \int_0^\infty \mathbf{P}(\xi\leq t)e^{-\theta t}\d t \\
& = 1-\theta\int_0^\infty\mathbf{P}(\xi>t)e^{-\theta t}\d t \\
& = 1-\theta\mathbf{E}[\xi]+\theta\int_0^\infty\mathbf{P}(\xi>t)(1-e^{-\theta t})\d t.
\end{align*}
Since $1-e^{-t}\leq \min\{1,t\}$,
\begin{align*}
\theta\int_0^\infty\mathbf{P}(\xi>t)(1-e^{-\theta t})\d t
& = \int_0^\infty\mathbf{P}(\xi>t/\theta)(1-e^{-t})\d t \\
& \leq  \int_0^\theta t\d t+\int_\theta^\infty\min\{1,t\}\mathbf{P}(\xi>t/\theta)\d t,
\end{align*}
and so setting $\theta=1/a_n$ and applying \eqref{e:upp} yields
\begin{align*}
\theta\int_0^\infty\mathbf{P}(\xi>t)(1-e^{-\theta t})\d t& \leq a_n^{-2}+ \frac{L(a_n)}{a_n^\alpha} \int_{a_n^{-1}}^\infty \min\{1,t\}t^{-\alpha}\frac{L(ta_n)}{L(a_n)}\d t.
\end{align*}
Now, by the definition of $a_n$, we have that $nL(a_n)a_n^{-\alpha}$ converges. Furthermore, by Theorem \ref{t:Pot}, we have that for any $\delta>0$, we can choose $A$ sufficiently large such that $L(ta_n)/L(a_n)$ is bounded above by $A\max\{t^\delta,t^{-\delta}\}$, uniformly in $n$ and $t\geq a_n^{-1}$. Since $\alpha\in(1,2)$, by taking $\delta$ suitably small so that $\delta-\alpha<-1<-\delta+1-\alpha$, we have that
\begin{align*}
\int_{a_n^{-1}}^\infty \min\{1,t\}t^{-\alpha}\frac{L(ta_n)}{L(a_n)}\d t\leq A\left(\int_0^1 t^{-\delta+1-\alpha}\d t +\int_1^\infty t^{\delta-\alpha}\d t\right)
\end{align*}
is bounded above by a finite constant uniformly in $n$. It also holds that $a_n^{-2}\leq Cn^{-1}$. Therefore,
\begin{equation*}
\theta\int_0^\infty\mathbf{P}(\xi>t)(1-e^{-\theta t})\d t \leq Cn^{-1}.
\end{equation*}
In particular, inserting this into the bounds above, we find that
\begin{align}
\mathbf{P}\left(\sum_{i=1}^n(\xi_i-\mathbf{E}[\xi_i])<-\lambda a_n\right)
& \leq \left(1-\frac{\mathbf{E}[\xi]}{a_n}+Cn^{-1}\right)^n\exp\left(n\frac{\mathbf{E}[\xi]}{a_n}-\lambda\right)\notag\\
& \leq \exp\left(-n\frac{\mathbf{E}[\xi]}{a_n}+C\right)\exp\left(n\frac{\mathbf{E}[\xi]}{a_n}-\lambda\right)\notag\\
& = \exp\left(C-\lambda\right), \label{e:neg}
\end{align}
which gives a bound that is better than polynomial for the probability of seeing lower fluctuations.

We next consider the upper fluctuations. Let $h$ be an increasing divergent positive function then, by \eqref{e:upp},
\begin{align}
\mathbf{P}\left(\bigcup_{i=1}^n\{\xi_i>h(\lambda)a_n\}\right)
& \leq n\mathbf{P}(\xi>h(\lambda)a_n) \notag\\
& \leq n(h(\lambda) a_n)^{-\alpha}L(h(\lambda)a_n),  \label{e:fUni}
\end{align}
which, by the choice of $a_n$, converges to $Ch(\lambda)^{-\alpha}$ as $n\rightarrow\infty$. Moreover, by Theorem \ref{t:Pot}, for any $\delta>0$, there exists $A$ such that \eqref{e:fUni} is bounded above by $Ah(\lambda)^{\delta-\alpha}$, uniformly in $n,\lambda\geq 1$. If we set, $\theta=(h(\lambda)a_n)^{-1}$, then it further holds that
\begin{align*}
\mathbf{E}\left[\exp\left(\theta \xi\right)\ind_{\xi\leq \theta^{-1}}\right]
& = \int_0^\infty \mathbf{P}\left(\exp\left(\theta \xi\right)\ind_{\xi\leq \theta^{-1}}>t\right)\d t \\
&\leq\theta\int_{-\infty}^{\theta^{-1}}\mathbf{P}(\xi>t)e^{\theta t}\d t \\
& \leq\theta\int_{-\infty}^0e^{\theta t}\d t+\theta\int_0^{\theta^{-1}}\mathbf{P}(\xi>t)\d t+\theta\int_0^{\theta^{-1}}\mathbf{P}(\xi>t)(e^{\theta t}-1)\d t \\
& \leq 1+\theta \mathbf{E}[\xi]+\theta\int_0^{\theta^{-1}}\mathbf{P}(\xi>t)(e^{\theta t}-1)\d t.
\end{align*}
Using \eqref{e:upp} and that $e^t\leq 1+ t + t^2\leq 1+2t$ for $t\in(0,1)$, we have that
\begin{align*}
\theta\int_0^{\theta^{-1}}\mathbf{P}(\xi>t)(e^{\theta t}-1)\d t
& \leq \theta\int_0^1(e^{\theta t}-1) \d t +\int_\theta^1\mathbf{P}(\xi>\theta^{-1} t)(e^t-1) \d t \\
& \leq \theta^2+2\theta^\alpha L(\theta^{-1}) \int_\theta^1t^{1-\alpha}\frac{L(\theta^{-1}t)}{L(\theta^{-1})}\d t.
\end{align*}
By Theorem \ref{t:Pot} we have that the above integral is bounded above uniformly in $n,\lambda\geq 1$, and, again by the choice of $a_n$, we have that $n\theta^\alpha L(\theta^{-1})$ converges to $Ch(\lambda)^{-\alpha}$ as $n\rightarrow \infty$. Moreover, by Theorem \ref{t:Pot}, for any $\delta>0$ there exists $A$ such that $n\theta^\alpha L(\theta^{-1})$ is bounded above by $Ah(\lambda)^{\delta-\alpha}$, uniformly in $n,\lambda\geq1$. In particular, we have that
\begin{equation*}
\theta\int_0^{\theta^{-1}}\mathbf{P}(\xi>t)(e^{\theta t}-1)\d t \leq \theta^2+Cn^{-1}h(\lambda)^{\delta-\alpha} \leq Cn^{-1}h(\lambda)^{\delta-\alpha}
\end{equation*}
for all $n,\lambda\geq 1$. Therefore, by Markov's inequality and \eqref{e:fUni},
\begin{align*}
&\mathbf{P}\left(\sum_{i=1}^n(\xi_i-\mathbf{E}[\xi_i])>\lambda a_n\right)  \\
&\quad \leq \mathbf{P}\left(\bigcup_{i=1}^n\{\xi_i>h(\lambda)a_n\}\right)+\mathbf{P}\left(\sum_{i=1}^n(\xi_i-\mathbf{E}[\xi_i])>\lambda a_n,\bigcap_{i=1}^n\{\xi_i\leq \theta^{-1}\}\right)  \\
& \quad \leq Ch(\lambda)^{\delta-\alpha}+ \mathbf{E}\left[\exp\left(\theta \xi\right)\ind_{\xi\leq \theta^{-1}}\right]^n\exp\left(-n\theta\mathbf{E}[\xi]-\theta\lambda a_n\right) \\
& \quad \leq Ch(\lambda)^{\delta-\alpha}+\left(1+\frac{\mathbf{E}[\xi]}{h(\lambda)a_n}+Cn^{-1}h(\lambda)^{\delta-\alpha}\right)^n\exp\left(-n\frac{\mathbf{E}[\xi]}{h(\lambda)a_n}-\lambda/h(\lambda)\right) \\
& \quad  \leq Ch(\lambda)^{\delta-\alpha}+\exp\left(Ch(\lambda)^{\delta-\alpha}-\lambda/h(\lambda)\right).
\end{align*}
Choosing $h(\lambda)=\lambda^\eta$ for some $\eta\in(0,1)$ (and $\delta$ suitably small), this yields that
\[\mathbf{P}\left(\sum_{i=1}^n(\xi_i-\mathbf{E}[\xi_i])>\lambda a_n\right) \leq C(\lambda^{-(\alpha-\delta)\eta} + e^{-\lambda^{1-\eta}}),\]
uniformly in $n,\lambda\geq 1$.

Combining the conclusion of the previous paragraph with \eqref{e:neg}, for any $\rho< \alpha$, we can choose $C_\rho$ such that
\begin{equation}\label{e:utp}
\mathbf{P}\left(\left|\sum_{i=1}^n\frac{\xi_i-\mathbf{E}[\xi_i]}{a_n}\right|>\lambda\right)\leq C_\rho\lambda^{-\rho}
\end{equation}
uniformly in $n,\lambda\geq 1$ for non-negative random variables $\xi$.

For general $\Rb$-valued random variables $\xi$ satisfying \eqref{e:upp} write $\xi^+=\xi\ind_{\xi\geq 0}$ and $\xi^-=|\xi|\ind_{\xi\leq 0}$. Using the triangle inequality and \eqref{e:utp} we have that
\begin{align*}
\lefteqn{\mathbf{P}\left(\left|\sum_{i=1}^n\frac{\xi_i-\mathbf{E}[\xi_i]}{a_n}\right|>\lambda\right)}\\
& = \mathbf{P}\left(\left|\sum_{i=1}^n\frac{\xi_i^+-\mathbf{E}[\xi_i^+]}{a_n}-\sum_{i=1}^n\frac{\xi_i^--\mathbf{E}[\xi_i^-]}{a_n}\right|>\lambda\right) \\
& \leq \mathbf{P}\left(\left|\sum_{i=1}^n\frac{\xi_i^+-\mathbf{E}[\xi_i^+]}{a_n}\right|>\lambda/2\right)+\mathbf{P}\left(\left|\sum_{i=1}^n\frac{\xi_i^--\mathbf{E}[\xi_i^-]}{a_n}\right|>\lambda/2\right)\\
& \leq 2^{1+\rho}C_\rho\lambda^{-\rho}
\end{align*}
uniformly in $n,\lambda\geq 1$. The tightness result readily follows.
\end{proof}

\begin{proof}[Proof of Proposition \ref{p:Lower}]
By Lemma \ref{l:nond}, we may assume that $t^{\alpha}/L(t)$ is non-decreasing. Let $(\xi'_i)_{i\geq 1}$ be a sequence of i.i.d.\ random variables which are independent of $(\xi_i)_{i\geq 1}$ and with the same law, and write $S'_n=\sum_{i=1}^n\xi'_i$ for their partial sum. We then have that
\begin{align*}
\mathbf{P}\left(\left|\frac{S_n-c_n}{a_n}\right|<\lambda\right)
& = \left(\mathbf{P}\left(\left|\frac{S_n-c_n}{a_n}\right|<\lambda\right)^2\right)^{1/2} \\
& = \mathbf{P}\left(\left|\frac{S_n-c_n}{a_n}\right|<\lambda,\left|\frac{S'_n-c_n}{a_n}\right|<\lambda \right)^{1/2} \\
& \leq \mathbf{P}\left(\left|\frac{S_n-S'_n}{a_n}\right|<2\lambda\right)^{1/2}.
\end{align*}
Writing $\overline{\xi}_i=\xi_i-\xi'_i$ and $\overline{S}_n=S_n-S'_n$, it thus suffices to show that
\begin{equation}\label{e:bound}
\mathbf{P}\left(\left|\overline{S}_n\right|<\lambda a_n\right)\leq C_\rho \lambda^{2\rho}
\end{equation}
for $n\geq 1$ and $\lambda\geq 1/a_n$.

Since $\mathbf{P}(|\xi|<t)\geq 1/2$ for $t$ sufficiently large, by \eqref{e:low} we have that
\begin{align*}
\mathbf{P}(|\overline{\xi}|>t)
& \geq \mathbf{P}(|\xi|>2t,|\xi'|<t) \\
& = \mathbf{P}(|\xi|>2t)\mathbf{P}(|\xi'|<t) \\
& \geq 2^{-(\alpha+1)}t^{-\alpha}L(2t) \\
& \geq Ct^{-\alpha}L(t)
\end{align*}
for $t$ sufficiently large. In particular, we can choose a constant $C$ such that $\mathbf{P}(|\overline{\xi}|>t) \geq Ct^{-\alpha}L(t)$ for all $t\geq 1$. By \cite[Theorem III.4]{P}, we have that there exists a universal constant $C<\infty$ such that
\[\sup_x\mathbf{P}(\overline{S}_n\in(x,x+t))\leq \frac{C}{\sqrt{n\mathbf{P}(|\overline{\xi}|>t/2)}}\]
for any $t\geq 0$. Choosing $t=2\lambda a_n$ gives
\[\mathbf{P}(|\overline{S}_n|\leq \lambda a_n)\leq \frac{C}{\sqrt{n\mathbf{P}(|\overline{\xi}|>\lambda a_n)}}\leq C\lambda^{\alpha/2}\left(\frac{a_n^\alpha}{nL(a_n)}\cdot \frac{L(a_n)}{L(\lambda a_n)}\right)^{1/2}\]
for $n\geq 1$ and $\lambda\geq 1/a_n$. Now, we have that $a_n^\alpha/nL(a_n)$ converges as $n\rightarrow\infty$ by the choice of $a_n$. Since $1/L(n)$ is locally bounded on $[1,\infty)$, by Theorem \ref{t:Pot}, we moreover have that for any $\delta>0$ there exists $A_\delta$ such that $L(a_n)/L(\lambda a_n)\leq A_\delta \lambda^{-\delta}$ for $\lambda\in[a_n^{-1},1]$. Taking $\delta$ small, the bound at \eqref{e:bound} follows.

Note that since $1/a_n\rightarrow0$, the given probability bound implies that
\[\lim_{\lambda\rightarrow0}\limsup_{n\rightarrow\infty}\mathbf{P}(|\overline{S}_n|\leq \lambda a_n)=0,\]
which completes the tightness result.
\end{proof}

\begin{proof}[Proof of Corollary \ref{c:convProb}]
If $\gamma\in(1,2)$ and \eqref{upperandlower} holds, then, by Propositions \ref{p:Upper} and \ref{p:Lower}, for any $\varepsilon>0$, there exists $\rho>0$ and a constant $C_{\rho,\delta}$ such that, for all $n\geq1$,
\begin{align}
\mathbf{P}\left(\left|\frac{ S_n-n\mathbf{E}\xi}{n^{\frac{1}{(1-\varepsilon)\gamma}}}\right|\geq \lambda\right)&\leq C_{\rho,\delta}\lambda^{-\rho},\qquad \forall \lambda\geq 1,\label{e:LowBnd} \\
\mathbf{P}\left(\left|\frac{ S_n-n\mathbf{E}\xi}{n^{\frac{1}{(1+\varepsilon)\gamma}}}\right|\leq \lambda\right)&\leq C_{\rho,\delta}\lambda^{\rho},\qquad \forall \lambda\geq n^{-1/(1+\varepsilon)\gamma}.\label{e:UppBnd}
\end{align}
The desired conclusion readily follows.
\end{proof}

\subsection{Almost-sure asymptotics}

In this subsection, we establish that, under the same conditions as Corollary \ref{c:convProb}, one does not have a corresponding almost-sure limit. Indeed, whilst the almost-sure limsup of the centred and rescaled partial sums matches the limit seen in  Corollary \ref{c:convProb}, almost-sure fluctuations mean that the almost-sure liminf is zero.

\begin{prop}\label{p:limsup}
Let $\gamma\in(1,2)$. Suppose that \eqref{upperandlower} holds. It is then the case that, $\mathbf{P}$-a.s.,
\[\limsup_{n\rightarrow\infty}\frac{\log \left| S_n-n\mathbf{E}\xi\right|}{\log n}=\frac{1}{\gamma},\]
\[\liminf_{n\rightarrow\infty}\frac{\log_+ \left| S_n-n\mathbf{E}\xi\right|}{\log n}=0.\]
\end{prop}
\begin{proof}
To establish the limsup result, we will show that, for any $\varepsilon>0$,
\begin{align}
\mathbf{P}\left(\limsup_{n\rightarrow\infty}\frac{\log \left| S_n-n\mathbf{E}\xi\right|}{\log n}<\frac{1}{\gamma}-\varepsilon\right)&=0\label{e:LowLog} \\
\mathbf{P}\left(\limsup_{n\rightarrow\infty}\frac{\log \left| S_n-n\mathbf{E}\xi\right|}{\log n}> \frac{1}{\gamma}+\varepsilon\right)&=0\label{e:UppLog}.
\end{align}
Clearly it will be sufficient to consider the case when $\varepsilon<1/\gamma$, which means that $n^{-1/\gamma}=o(n^{-\varepsilon})$. For \eqref{e:LowLog}, we have, for any $\delta>0$, that
\begin{align*}
\mathbf{P}\left(\limsup_{n\rightarrow\infty}\frac{\log \left| S_n-n\mathbf{E}\xi\right|}{\log n}< \frac{1}{\gamma}-\varepsilon\right)
& \leq \limsup_{n\rightarrow\infty}\mathbf{P}\left(\left|\frac{ S_n-n\mathbf{E}\xi}{n^{\frac{1}{(1+\delta)\gamma}}}\right|\leq n^{\frac{\delta-\varepsilon(1+\delta)\gamma}{(1+\delta)\gamma}}\right).
\end{align*}
On choosing $\delta$ sufficiently small, the result thus follows from \eqref{e:UppBnd}. For \eqref{e:UppLog}, let $b>1$ and define $n_k=b^k$. It is then the case that
\begin{align*}
&\mathbf{P}\left(\limsup_{n\rightarrow\infty}\frac{\log \left| S_n-n\mathbf{E}\xi\right|}{\log n}> \frac{1}{\gamma}+\varepsilon\right)\\
& \leq\mathbf{P}\left(\limsup_{n\rightarrow\infty}\frac{\left| S_n-n\mathbf{E}\xi\right|}{n^{\frac{1}{\gamma}+\varepsilon}}\geq 1\right) \\
& \leq \mathbf{P}\left(\limsup_{k\rightarrow\infty}\frac{\left| S_{n_k}-n_k\mathbf{E}\xi\right|}{n_k^{\frac{1}{\gamma}+\varepsilon}}\geq 1/2\right)\\
&+  \mathbf{P}\left(\bigcap_{m\geq1}\bigcup_{k\geq m}\sup_{l\in\{n_k,...,n_{k+1}-1\}}\left|\frac{\left| S_{l}-l\mathbf{E}\xi\right|}{l^{\frac{1}{\gamma}+\varepsilon}}-\frac{\left| S_{n_k}-n_k\mathbf{E}\xi\right|}{n_k^{\frac{1}{\gamma}+\varepsilon}}\right|\geq 1/2\right).
\end{align*}
By \eqref{e:LowBnd} and choosing $\delta$ sufficiently small so that $\varepsilon(1-\delta)\gamma-\delta>0$, we have that there exists $\eta>0$ such that
\begin{align*}
\mathbf{P}\left(\frac{\left| S_{n_k}-n_k\mathbf{E}\xi\right|}{n_k^{\frac{1}{\gamma}+\varepsilon}}\geq 1/2\right)
 = \mathbf{P}\left(\frac{\left| S_{n_k}-n_k\mathbf{E}\xi\right|}{n_k^{\frac{1}{(1-\delta)\gamma}}}\geq\frac{n_k^{\frac{1}{(1-\delta)\gamma}(\varepsilon(1-\delta)\gamma-\delta)}}{2}\right)
 \leq C_\eta b^{-\eta k}.
\end{align*}
By the Borel-Cantelli lemma, we therefore have that
\begin{equation}\label{bbb}
\mathbf{P}\left(\limsup_{k\rightarrow\infty}\frac{\left| S_{n_k}-n_k\mathbf{E}\xi\right|}{n_k^{\frac{1}{\gamma}+\varepsilon}}\geq 1/2\right)=0.
\end{equation}
For $l\in\{n_k,...,n_{k+1}-1\}$ we have that
\[\left|\frac{\left| S_{l}-l\mathbf{E}\xi\right|}{l^{\frac{1}{\gamma}+\varepsilon}}-\frac{\left| S_{n_k}-n_k\mathbf{E}\xi\right|}{n_k^{\frac{1}{\gamma}+\varepsilon}}\right|
\leq \frac{\left|S_{l}-S_{n_k}-(l-n_k)\mathbf{E}\xi\right|}{n_k^{\frac{1}{\gamma}+\varepsilon}}
+\left| S_{n_k}-n_k\mathbf{E}\xi\right|\left(\frac{1}{n_k^{\frac{1}{\gamma}+\varepsilon}}-\frac{1}{n_{k+1}^{\frac{1}{\gamma}+\varepsilon}}\right), \]
where we note that
\[\frac{1}{n_k^{\frac{1}{\gamma}+\varepsilon}}-\frac{1}{n_{k+1}^{\frac{1}{\gamma}+\varepsilon}}=\frac{b^{\frac{1}{\gamma}+\varepsilon}-1}{b^{\frac{1}{\gamma}+\varepsilon}}\cdot\frac{1}{n_k^{\frac{1}{\gamma}+\varepsilon}}\leq \frac{1}{2n_k^{\frac{1}{\gamma}+\varepsilon}} \]
for $b>1$ sufficiently small. In particular, applying this bound with \eqref{bbb}, we deduce that
\[ \mathbf{P}\left(\limsup_{k\rightarrow\infty}\left| S_{n_k}-n_k\mathbf{E}\xi\right|\left(\frac{1}{n_k^{\frac{1}{\gamma}+\varepsilon}}-\frac{1}{n_{k+1}^{\frac{1}{\gamma}+\varepsilon}}\right)\geq 1/4\right)=0.\]
Consequently, it remains to show that
\begin{equation}\label{aaa}
 \mathbf{P}\left(\bigcap_{m\geq1}\bigcup_{k\geq m}\sup_{l\in\{n_k,...,n_{k+1}-1\}}\frac{\left|S_{l}-S_{n_k}-(l-n_k)\mathbf{E}\xi\right|}{n_k^{\frac{1}{\gamma}+\varepsilon}}\geq 1/4\right)=0.
 \end{equation}
Note that $S_l-l\mathbf{E}\xi$ is a martingale, and applying Doob's inequality to it yields
\begin{align*}
&\mathbf{P}\left(\sup_{l\in\{n_k,...,n_{k+1}-1\}}\frac{\left|S_{l}-S_{n_k}-(l-n_k)\mathbf{E}\xi\right|}{n_k^{\frac{1}{\gamma}+\varepsilon}}\geq 1/4\right)\\
& \leq \mathbf{P}\left(\sup_{l\in\{1,...,(b-1)n_{k}\}}\frac{\left|S_{l}-l\mathbf{E}\xi\right|}{n_k^{\frac{1}{\gamma}+\varepsilon}}\geq 1/4\right)\\
& \leq \frac{C_p\mathbf{E}\left|S_{(b-1)n_{k}}-(b-1)n_{k}\mathbf{E}\xi\right|^p}{n_k^{p\left(\frac{1}{\gamma}+\varepsilon\right)}}
\end{align*}
for any $p>1$. Now, by \eqref{upperandlower} we have that $\mathbf{E}|\xi|^p<\infty$ for all $p<\gamma$. Hence, since $\xi_k-\mathbf{E}\xi$ are centred, independent and identically distributed, by \cite[Theorem 2]{BE} we have that
\[\mathbf{E}\left|\sum_{k=1}^n(\xi_k-\mathbf{E}\xi)\right|^p\leq 2\sum_{k=1}^n\mathbf{E}|\xi_k-\mathbf{E}\xi|^p\leq 4n\mathbf{E}|\xi|^p\]
 for any such $p$. Choosing $p<\gamma$ such that $p(\frac1\gamma+\varepsilon)>1$, we then have that
\begin{align*}
 \frac{\mathbf{E}\left|S_{(b-1)n_{k}}-(b-1)n_{k}\mathbf{E}\xi\right|^p}{n_k^{p\left(\frac{1}{\gamma}+\varepsilon\right)}}
  \leq \frac{Cn_{k}\mathbf{E}|\xi|^p}{n_k^{p\left(\frac{1}{\gamma}+\varepsilon\right)}}
  \leq Cb^{-k\eta}
\end{align*}
for some $\eta>0$. The result at \eqref{aaa} then follows by the Borel-Cantelli lemma, and this is enough to complete the proof of the limsup part of the proposition.

We now turn to the liminf. Call a value $x\in\Rb$ recurrent if we have that
\[\mathbf{P}\left(\liminf_{n\rightarrow\infty}|S_n-n\mathbf{E}\xi-x|=0\right)=1.\]
By \eqref{upperandlower}, we have that $\mathbf{E}|\xi-\mathbf{E}\xi|<\infty$. Therefore, by \cite[Theorem 4]{CF}, we have that $0$ is recurrent, which proves that
\[\liminf_{n\rightarrow\infty}\frac{\log_+ \left| S_n-n\mathbf{E}\xi\right|}{\log n}=0.\]
\end{proof}

\section{Biased random walk on supercritical percolation}\label{spsec}

In this section, we study the model of biased random walk on supercritical percolation. We begin, in Section \ref{s:reg}, by describing the regeneration structure used in \cite{FH,S} and building on the estimates proved therein. In Section \ref{s:conc}, we then combine these estimates for regeneration times with the results of Section \ref{nssec} to prove Theorem \ref{mainthm}. We note that the only results with proofs that depend on the particular setting of biased random walk on supercritical percolation are those introduced in Section \ref{s:reg}. Specifically, for a given model of a ballistic random walk with a suitable regeneration structure satisfying statements analogous to Theorems \ref{t:Xexp}, \ref{t:up} and \ref{t:lower}, the conclusions of Theorem \ref{mainthm} will follow.

\subsection{Regeneration times}\label{s:reg}

Our proof of Theorem \ref{mainthm} will be based on the (configuration dependent) regeneration times of the random walk in random environment used in \cite{FH,S}. These times separate the trajectory of the walk into independent blocks; this allows us to exploit our results from the previous section to prove Theorem \ref{mainthm}. We now construct these formally following \cite{FH,S}.

Throughout, we assume that $p>p_c(d)$ so that there exists, $P_p$-a.s., a unique infinite open cluster and we write $\Ic_x$ for the event that $x$ is contained in this cluster. For $x\in\Zb^d$, $\omega\in \Ic_x$ and $j\geq 0$, let $\Pc_{j,x}(\omega)$ denote the set of simple paths $(\pi(i))_{i\geq 1}$ in the unique infinite open cluster $\Cc_\infty$ started from $\pi(0)=x$ such that $\pi(i)\in\Hc_x^+$ for $i\geq j$, where
\[\Hc_x^+=\{y\in\Zb^d: y\cdot\vec{\ell}\geq x\cdot\vec{\ell} \}.\]
Now, $P_p$-a.s., on $\Ic_x$, $\Pc_{j,x}(\omega)$ is non-empty for large enough $j$, and we define
\[J_x(\omega)=\begin{cases}\inf\{j\geq 0: \Pc_{j,x}(\omega)\neq \emptyset\}, & \text{if } \omega\in\Ic_x,\\ \infty,&\text{if }\omega\notin\Ic_x.\end{cases}\]
Note that this is similar to the backtracking function of \eqref{bt}, but based on the number of steps taken to escape the region lower than $x$ in the direction of $\vec{\ell}$, rather than the depth of the trap based at $x$. We then define $W_0=0$, $m_0=J_{X_0}(\omega)$ and, by induction,
\[\begin{cases}W_{k+1}=2+\Delta_{m_k}, \\ m_{k+1}=\sup\{X_n\cdot\vec{\ell}:n\leq W_{k+1}\}+1, \end{cases}\]
for all $k\geq 0$. Next, let $(e_i)_{i=1}^d$ be an orthonormal basis of $\Zb^d$ such that $e_1\cdot\vec{\ell}\geq e_2\cdot\vec{\ell}\geq\ldots\geq e_d\cdot\vec{\ell}\geq 0$, and define $B$ to be the collection of edges of the form $[-e_1,e-e_1]$ where $e$ is any unit vector satisfying $e\cdot\vec{\ell}=e_1\cdot\vec{\ell}$. We then define the stopping times
\begin{align*}
\sigma_1&=\inf\{W_k:k\geq 1, X_{W_k}=X_{W_k-1}+e_1=X_{W_k-2}+2e_1, \omega(b)=1,  \forall b\in B+X_{W_k}\},\\
D&=\inf\{n\geq0:X_n\cdot\vec{\ell}<X_0\cdot\vec{\ell}\}.
\end{align*}
Writing $M_0=X_0\cdot\vec{\ell}$, we recursively define
\[\begin{cases}
\sigma_{k+1}=\sigma_1\circ \theta_{\Delta_{M_k}}, \\
R_{k+1}=D\circ \theta_{\sigma_{k+1}}+\sigma_{k+1}, \\
M_{k+1}=\sup_{n\leq R_{k+1}}X_n\cdot\vec{\ell}+1
\end{cases}\]
for $k\geq 0$. Finally, we define the first regeneration time $\tau_1=\sigma_K$, where
\[K=\inf\{k\geq 1:\sigma_k<\infty \;\text{and}\;R_k=\infty\},\]
and the increasing sequence of regeneration times $(\tau_k)_{k\geq 1}$ recursively via
\[\tau_{k+1}=\tau_1+\tau_k(X_{\tau_1+\cdot}-X_{\tau_1}, \omega(\cdot+X_{\tau_1}))\quad \text{for } k\geq 1.\]

By \cite[Lemma 2.2 and Theorem 2.4]{S} (as quoted in \cite[Lemma 4.1 and Theorem 4.2]{FH}), the sequence $(\tau_k)_{k\geq 1}$ exists such that, for any $k$, $P_p$-a.s.\ for all $x\in \Cc_\infty(\omega)$, $\tau_k<\infty$ $P_x^\omega$-a.s.\ and, under $\Pb_p$, the processes
\[(X_{\min\{\tau_1, \cdot\}}),(X_{\min\{\tau_1+\cdot,\tau_2\}}-X_{\tau_1}),(X_{\min\{\tau_2+\cdot, \tau_3\}}-X_{\tau_2}),\ldots\]
are independent and, except for the first, identically distributed. In particular, it follows that $\tau_1, \tau_2-\tau_1, \tau_3-\tau_2,\ldots$ are independent and, except for the first, identically distributed.

In \cite{FH}, the following upper bound was given on the tail of the distance between regeneration points. This shows that the walk does not travel a large distance between regeneration times.

\begin{thm}[{\cite[Proposition 4.4]{FH}}]\label{t:Xexp}
Let $d\geq 2$, $p\in(p_c(d),1)$ and $\vec{\ell}\in\mathbb{S}^{d-1}$, there exist constants $C>c>0$ such that, for all $t>0$,
\begin{align*}
\Pb_p((X_{\tau_2}-X_{\tau_1})\cdot\vec{\ell}\geq t)\leq Ce^{-ct}.
\end{align*}
\end{thm}

The anomalous behaviour is driven by trapping, for which an upper bound is given by the following key estimate on the tail of the regeneration times of \cite{FH}. Applying Proposition \ref{p:Upper}, this is already enough to give an upper bound on the fluctuations of $\tau_n$ about its mean of the desired order.

\begin{thm}[{\cite[Theorem 4.3]{FH}}]\label{t:up} Let $d\geq 2$, $p\in(p_c(d),1)$ and $\vec{\ell}\in\mathbb{S}^{d-1}$. If $\gamma$ is defined as at \eqref{gammadef} and $\varepsilon\in(0,1)$, then there exists a constant $C\in(0,\infty)$ such that
\[\mathbb{P}_p\left(\tau_2-\tau_1\geq t\right)\leq Ct^{-(1-\varepsilon)\gamma},\qquad \forall t\geq 1.\]
\end{thm}

The first main step of this section is to prove the corresponding lower bound, as we now present.

\begin{thm}\label{t:lower} Let $d\geq 2$, $p\in(p_c(d),1)$ and $\vec{\ell}\in\mathbb{S}^{d-1}$. If $\gamma$ is defined as at \eqref{gammadef} and $\varepsilon\in(0,1)$, then there exists a constant $C\in(0,\infty)$ such that
\[\mathbb{P}_p\left(\tau_2-\tau_1\geq t\right)\geq Ct^{-(1+\varepsilon)\gamma},\qquad \forall t\geq 1.\]
\end{thm}

Our proof of Theorem \ref{t:lower} will be based on the excursion times in one-headed traps of \cite{FH}. For an environment $\omega$ and vertices $x\sim y$, write $\omega([x,y])$ for the environment $\omega$ except that the edge $[x,y]$ is closed. We moreover write $K^\omega(x)$ for the open cluster containing $x$. We say that there is a one-headed trap $\Tc(x)$ with head $x$ if
\begin{enumerate}
\item $[x,x+e_1]$ is open;
\item $|K^{\omega([x,x+e_1])}(x+e_1)|<\infty$;
\item $(x+e_1)\cdot\vec{\ell}\leq y\cdot\vec{\ell}$ for $y\in K^{\omega([x,x+e_1])}(x+e_1)$.
\end{enumerate}
If these hold, then we set $\Tc(x)=K^{\omega([x,x+e_1])}(x+e_1)$, otherwise $\Tc(x)=\emptyset$. (See \cite[Figure 3.1]{FH}.) We then write $\delta(x)$ for the vertex $y\in\Tc(x)$ maximising $y\cdot\vec{\ell}$ which is chosen according to some predetermined order on $\Zb^d$ if there are several such vertices. The furthest distance the walk can reach in direction $\vec{\ell}$ within a trap rooted at a vertex $x$ is
\[\Dk(x)=\begin{cases}0, & \text{if } \omega \in \Ic^+_x, \\ \max_{\Pc\in\Pc_x^+}\max_{y\in\Pc}(y-x)\cdot\vec{\ell}, & \text{otherwise},\end{cases}\]
where $\Pc_x^+$ is the set of simple paths starting from $x$ that remain in $\Hc_x^+=\{z\in\Zb^d:z\cdot\vec{\ell}\geq x\cdot\vec{\ell}\}$, and $\Ic^+_x$ is the subset of $\Omega$ consisting of configurations such that $x$ belongs to the infinite connected component induced by the restriction to edges in $\Hc_x^+$. We next subdivide space into \emph{slabs}, writing
\[\Upsilon_{k,n}=\{x\in\Zb^d:x\cdot\vec{\ell}\in[k\log(n)^3,(k+1)\log(n)^3) \}\]
for the $k$th slab at scale $n$. We set
\[T_{k,n}=\inf\{m\geq0:X_m\in \Upsilon_{k,n}\}\]
to be the time taken by the random walk to reach the $k$th slab, and let $Y_{k,n}=X_{T_{k,n}}$ be the entrance point of the $k$th slab. For $k,n\in\mathbb{N}$ and $\varepsilon>0$ fixed, we write
\[Z_{k,n}=\left\{\Tc(Y_{k,n})\neq \emptyset,\:\Dk(Y_{k,n}+e_1)\geq (1-\varepsilon)\zeta^{-1}\log n,\:|\Tc(Y_{k,n})|\leq (\log n)^2\right\}\]
for the event that the trap $\Tc(Y_{k,n})$ is suitably deep (and not too large). Furthermore, let
\[\tilde{Z}_{k,n}=Z_{k,n}\cap\{\inf\{m>T_{k,n}:X_m=Y_{k,n}\}>\inf\{m>T_{k,n}:X_m=\delta(Y_{k,n})\}\}\]
be the event that the walk \emph{explores} a deep one-headed trap upon entering the $k$th slab. Note that on this event it must be the case that $X_{T_{k,n}+1}=Y_{k,n}+e_1$ since the edge $[Y_{k,n},Y_{k,n}+e_1]$ is the only entrance into the trap. Finally, we write
\[\theta_{k,n}=\ind_{\tilde{Z}_{k,n}}(\inf\{m>T_{k,n}:X_m=Y_{k,n}\}-T_{k,n})\})\]
for the duration of the first excursion in the trap rooted at the entrance to the $k$th slab when it is a deep trap and explored by the walk. As described in \cite[Section 3]{FH}, for each $n$, $(\theta_{k,n})_{k=1}^{n/\log(n)^3}$ are independent lower bounds on the times taken for the walk to traverse the corresponding slabs. The proof of Theorem \ref{t:lower} will depend on bounds from \cite[Section 3]{FH} on the number of excursions into deep one-headed traps and the durations of these excursions.

\begin{proof}[Proof of Theorem \ref{t:lower}]
Since there can be at most $n$ regenerations up to distance $n$ in direction $\vec{\ell}$, and the time spent in deep one-headed traps before reaching distance $n$ in direction $\vec{\ell}$ is at most the time taken to reach distance $n$ in direction $\vec{\ell}$, it is immediate that
\[\sum_{k=1}^{n/\log(n)^3}\theta_{k,n}\leq \Delta_n\leq \tau_1+\sum_{k=1}^n(\tau_{k+1}-\tau_k).\]
For $M>1/c$, where $c$ is the constant of Theorem \ref{t:Xexp}, a union bound yields that
\begin{align*}
\Pb_p\left(\bigcup_{j=1}^n\{(X_{\tau_{j+1}}-X_{\tau_j})\cdot\vec{\ell}\geq M\log(n)\} \right)
\leq Cne^{-cM\log(n)},
\end{align*}
where the right-hand side converges to $0$ as $n\to\infty$. In particular, with high probability, none of the first $n$ `regeneration blocks' have width greater than $M\log(n)$. Since each deep one-headed trap is of depth at least $(1-\varepsilon)\zeta^{-1}\log(n)$ and the walk must reach the deepest point so that the trap is explored, it follows that, with high probability, there can be at most $M\zeta/(1-\varepsilon)$ deep one-headed traps explored in any of the first $n$ regeneration blocks. Note that we have used here that the walk can not regenerate during an excursion in a one-headed trap because the walk must return to the unique entrance to the trap in order to escape and, by definition, every vertex in the trap must be deeper in direction $\vec{\ell}$ than the entrance.

By \cite[Section 3]{FH}, with high probability, the walk encounters at least $n^{\varepsilon/2}$ deep one-headed traps in the first $n/\log (n)^3$ slabs. Moreover, conditional on encountering such a trap, it has a probability of at least $c/\log(n)^2$ of spending time at least $n^{(1-2\varepsilon)/\gamma}$ there (see \cite[Equation (3.7)]{FH} for a precise statement). In particular, the latter two observations mean that, with high probability, in at least $cn^{\varepsilon/2}/\log(n)^2$ of the first $n/\log (n)^3$ slabs, the walk spends time at least $n^{(1-2\varepsilon)/\gamma}$ in a deep one-headed trap. In conjunction with the conclusion of the previous paragraph, we thus can further conclude that, with high probability, in at least $cn^{\varepsilon/2}/\log(n)^2$ of the first $n$ regeneration blocks, the walk spends time at least $n^{(1-2\varepsilon)/\gamma}$. It follows that
\[\Pb_p(\tau_2-\tau_1\geq n^{(1-2\varepsilon)/\gamma})=n^{-1}\mathbb{E}_p\left(\sum_{k=1}^{n}\mathbf{1}_{\tau_{k+1}-\tau_{k}\geq n^{(1-2\varepsilon)/\gamma}}\right)\geq c\frac{n^{\varepsilon/2-1}}{\log(n)^2}\]
for $n$ suitably large. Using monotonicity, it is then straightforward to show that
\[\Pb_p(\tau_2-\tau_1\geq t)\geq ct^{-(1-\tilde{\varepsilon})\gamma}\]
for all $t\geq 1$ and some suitable $c,\tilde{\varepsilon}$.
\end{proof}

Putting together the results of Section \ref{nssec} and Theorems \ref{t:up} and \ref{t:lower}, we immediately obtain the following consequence concerning the fluctuations of regeneration times around their mean behaviour.

\begin{cor}\label{c:3conv} Let $d\geq 2$, $p\in(p_c(d),1)$ and $\vec{\ell}\in\mathbb{S}^{d-1}$. If $\gamma$ is defined as at \eqref{gammadef}, then
\begin{align*}
\frac{\log \left| \tau_{n+1}-\tau_1-n\mathbb{E}_p(\tau_2-\tau_1)\right|}{\log n}
&\buildrel{\mathbb{P}_p}\over\rightarrow\frac{1}{\gamma}, \\
\limsup_{n\to\infty}\frac{\log \left| \tau_{n+1}-\tau_1-n\mathbb{E}_p(\tau_2-\tau_1)\right|}{\log n}
& = \frac{1}{\gamma}, \qquad \Pb_p\text{-a.s.},\\
\liminf_{n\to\infty}\frac{\log_+ \left| \tau_{n+1}-\tau_1-n\mathbb{E}_p(\tau_2-\tau_1)\right|}{\log n}
& = 0 ,\qquad \Pb_p\text{-a.s.}
\end{align*}
\end{cor}

We are almost ready to prove Theorem \ref{mainthm}, but as a further ingredient for this, it will be helpful to derive a technical lemma that controls the fluctuations of the number of regenerations up to a given level. For $n\in\mathbb{N}$, write
\[\vartheta_n=\sup\{k\in\Nb: \tau_k\leq \Delta_n\}=\sup\{k\in\Nb:X_{\tau_k}\cdot\vec{\ell}\leq n\}\]
for the number of regenerations before reaching distance $n$ in direction $\vec{\ell}$. We then have that $(\vartheta_n)_{n\geq 1}$ is a non-decreasing sequence of $\Nb$-valued random variables such that $\vartheta_n-\vartheta_{n-1}\leq 1$ and
\begin{equation}\label{e:DeltaBound}
\tau_1+\sum_{k=1}^{\vartheta_n-1}(\tau_{k+1}-\tau_k)=\tau_{\vartheta_n}\leq \Delta_n< \tau_{\vartheta_n+1}=\tau_1+\sum_{k=1}^{\vartheta_n}(\tau_{k+1}-\tau_k).
\end{equation}
By Theorem \ref{t:Xexp} and the strong law of large numbers, we have that $X_{\tau_n}\cdot\vec{\ell}/n$ converges $\Pb_p$-a.s.\ to $\Eb_p[(X_{\tau_2}-X_{\tau_1})\cdot\vec{\ell}]$. It readily follows that $\vartheta_n/n$ converges $\Pb_p$-a.s.\ to $\eta=1/\Eb_p[(X_{\tau_2}-X_{\tau_1})\cdot\vec{\ell}]$. The next result, which is also a consequence of Theorem \ref{t:Xexp}, gives an estimate for the deviation between $\vartheta_n$ and $\eta n$. We use the standard abbreviation i.o.\ for `infinitely often'.

\begin{lem}\label{l:thetaDiff}
Let $d\geq 2$, $p\in(p_c(d),1)$ and $\vec{\ell}\in\mathbb{S}^{d-1}$. If $\gamma$ is defined as at \eqref{gammadef}, then for $D<\infty$ suitably large
\[\Pb_p(|\vartheta_n-n\eta|>Dn^{1/2}\log(n)\text{ i.o.})=0.\]
\end{lem}
\begin{proof}
Write $M_n=Dn^{1/2}\log(n)$, then
\[\{|\vartheta_n-n\eta|>M_n\}\subseteq\left\{\sum_{k=0}^{\lceil n\eta-M_n\rceil}(X_{\tau_{k+1}}-X_{\tau_k})\cdot\vec{\ell} < n < \sum_{k=0}^{\lfloor n\eta+M_n\rfloor}(X_{\tau_{k+1}}-X_{\tau_k})\cdot\vec{\ell} \right\}^c,\]
where for convenience we set $\tau_0=0$. Let $Y_k=(X_{\tau_{k+1}}-X_{\tau_k})\cdot\vec{\ell}-\Eb[(X_{\tau_{k+1}}-X_{\tau_k})\cdot\vec{\ell}]$ which, by Theorem \ref{t:Xexp}, has exponential moments for $k\geq 1$. Let $\tilde{M}_n=\eta^{-1}M_n-\log (n)$, then
\begin{align*}
\Pb_p\left(\sum_{k=1}^{\lceil n\eta-M_n\rceil}(X_{\tau_{k+1}}-X_{\tau_k})\cdot\vec{\ell}\geq n-\log(n)\right)
& \leq \Pb_p\left(\sum_{k=1}^{\lceil n\eta-M_n\rceil}Y_k\geq \tilde{M}_n\right) \\
& \leq  \Eb[e^{\lambda Y_1}]^{\lceil n\eta-M_n\rceil}e^{-\lambda \tilde{M}_n},
\end{align*}
where the upper bound is finite, for $\lambda>0$ suitably small. In particular, using a Taylor expansion and choosing $\lambda=\sqrt{2/n\eta}$, we have that the final term is bounded above by $Ce^{-c\tilde{M}_n/\sqrt{n}}\leq Ce^{-cD\log(n)}$. Choosing $D$ suitably large, we then have that
\[\sum_{n=1}^\infty \Pb_p\left(\sum_{k=1}^{\lceil n\eta-M_n\rceil}(X_{\tau_{k+1}}-X_{\tau_k})\cdot\vec{\ell} \geq n-\log (n)\right) <\infty\]
and so, by the Borel-Cantelli lemma and the fact that $X_{\tau_1}\cdot\vec{\ell}$ is almost-surely finite,
\[\Pb_p\left(\sum_{k=0}^{\lceil n\eta-M_n\rceil}(X_{\tau_{k+1}}-X_{\tau_k})\cdot\vec{\ell}\geq n\text{ i.o.}\right)=0.\]
Similarly, writing $\bar{M}_n=M_n\eta^{-1}$, for $\lambda>0$,
\begin{align*}
\Pb_p\left(\sum_{k=0}^{\lfloor n\eta+M_n\rfloor}(X_{\tau_{k+1}}-X_{\tau_k})\cdot\vec{\ell} \leq n\right)
& \leq \Pb_p\left(\sum_{k=1}^{\lfloor n\eta+M_n\rfloor}Y_k \leq-\bar{M}_n\right)\\
&\leq \Eb_p[e^{-\lambda Y_1}]^{\lfloor n\eta+M_n\rfloor}e^{-\lambda \bar{M}_n} .
\end{align*}
Again, choosing $\lambda=\sqrt{2/n\eta}$, we have that this is bounded above by $Ce^{-c\tilde{M}_n/\sqrt{n}}\leq Ce^{-cD\log(n)}$
and so, by the same argument,
\[\Pb_p\left(\sum_{k=0}^{\lfloor n\eta+M_n\rfloor}(X_{\tau_{k+1}}-X_{\tau_k})\cdot\vec{\ell} \leq n\text{ i.o.}\right)=0.\]
\end{proof}

\subsection{Proof of Theorem \ref{mainthm}}\label{s:conc}

In this subsection, we prove the claims of Theorem \ref{mainthm}. To this end, we will apply the results of Section \ref{nssec} to two sequences of random variables. One of these will be the sequence $(\tau_{m+1}-\tau_m-\Eb[\tau_2-\tau_1])_{m\geq 1}$, which was previously considered in Corollary \ref{c:3conv}. We further note that, by the strong law of large numbers, $\tau_{n+1}/n$ converges $\Pb_p$-a.s.\ to $\Eb_p[\tau_2-\tau_1]$, and therefore
\begin{align}\label{e:Ltn}
\lim_{n\to\infty}\frac{\log \tau_{n+1}}{\log n}=1, \qquad \Pb_p\text{-a.s.}
\end{align}
A second sequence that will appear in our proof of Theorem \ref{mainthm} is given by $(Z_m)_{m\geq 1}$, where
\[Z_m=X_{\tau_{m+1}}\cdot\vec{\ell}-X_{\tau_{m}}\cdot\vec{\ell}-v\left(\tau_{m+1}-\tau_{m}\right).\]
(Cf. \cite[Equation (4.6)]{S1}.) This is an i.i.d.\ sequence of mean 0 random variables. Moreover, by Theorems \ref{t:Xexp}, \ref{t:up} and \ref{t:lower}, we have that \eqref{upperandlower} holds for $Z_m$; that is, for every $\varepsilon\in(0,1)$, there exists $C\in(0,\infty)$ such that
\begin{align}\label{e:Zuplo}
C^{-1}t^{-(1+\varepsilon)\gamma}\leq \Pb_p(|Z_1|\geq t)\leq Ct^{-(1-\varepsilon)\gamma}, \qquad \forall t\geq 1.
\end{align}

We start with the proof of \eqref{deltalim}, but postpone the proof of the other limit in probability, i.e.\ \eqref{xlim}, until after we have checked the almost-sure parts, as it will be convenient to apply \eqref{limsup} in the proof of \eqref{xlim}.
\begin{proof}[Proof of \eqref{deltalim}]
Using \eqref{e:DeltaBound} and Markov's inequality, we have that
\begin{align*}
&\Pb_p(|\Delta_n-\tau_{\eta n}|>n^{\frac{1}{\gamma}-\varepsilon})\\
&  \leq \Pb_p(|\vartheta_n-n\eta|>Dn^{1/2}\log(n))
+\Pb_p\left(\sum_{k=\lfloor\eta n-Dn^{1/2}\log(n)\rfloor}^{\lceil\eta n+Dn^{1/2}\log(n)\rceil+1}(\tau_{k+1}-\tau_k)>n^{\frac{1}{\gamma}-\varepsilon}\right) \\
&  \leq \Pb_p(|\vartheta_n-n\eta|>Dn^{1/2}\log(n))
+\frac{2\lceil D n^{1/2}\log(n)\rceil+2}{n^{\frac{1}{\gamma}-\varepsilon}}\Eb_p[\tau_2-\tau_1].
\end{align*}
By Lemma \ref{l:thetaDiff}, for $\varepsilon<1/\gamma-1/2$, this converges to $0$ as $n\to\infty$. It then follows from the first statement of Corollary \ref{c:3conv} that
\begin{align*}
\frac{\log \left| \Delta_n-nv^{-1}\right|}{\log n}
&\buildrel{\mathbb{P}_p}\over\rightarrow\frac{1}{\gamma},
\end{align*}
and the result is readily extended to the continuous parameter $s$.
\end{proof}

\begin{proof}[Proof of \eqref{limsup}]
Let $\Sc_n=\sum_{m=1}^nZ_m$, then it follows from \eqref{e:Zuplo} and Proposition \ref{p:limsup} that, $\Pb_p$-a.s.,
\begin{equation}\label{e:lsZ}
\limsup_{n\to\infty}\frac{\log|\Sc_n|}{\log n}=\frac{1}{\gamma}
\end{equation}
Note that $\Sc_n=(X_{\tau_{n+1}}-X_{\tau_1})\cdot\vec{\ell}-v(\tau_{n+1}-\tau_1)$, where $X_{\tau_1}\cdot\vec{\ell}$ and $\tau_1$ are $\Pb_p$-a.s.\ finite and do not depend on $n$. By \eqref{e:Ltn}, we thus have that
\begin{align*}
\limsup_{n\to\infty}\frac{\log |X_{\tau_{n+1}}\cdot\vec{\ell}-v\tau_{n+1}|}{\log \tau_{n+1}}=\frac{1}{\gamma}, \qquad \Pb_p\text{-a.s.}
\end{align*}
and therefore, using \eqref{e:Ltn},
\begin{align*}
\limsup_{n\to\infty}\frac{\log |X_{n}\cdot\vec{\ell}-vn|}{\log n}\geq\frac{1}{\gamma}, \qquad \Pb_p\text{-a.s.}
\end{align*}

Write $\kappa_n=\sup\{m\geq 0:\tau_k\leq n\}$ for the number of regenerations by time $n$. By the law of large numbers for $\tau_n$, we have that $\kappa_n/n$ converges $\Pb_p$-a.s.\ to $1/\Eb_p[\tau_2-\tau_1]$, and therefore $\log \kappa_n/\log n$ converges $\Pb_p$-a.s.\ to $1$. It follows from \eqref{e:lsZ} that
\begin{equation}\label{e:lskappa}
\limsup_{n\to\infty}\frac{\log|\Sc_{\kappa_n}|}{\log n}=\limsup_{n\to\infty}\frac{\log|\Sc_{\kappa_n}|}{\log\kappa_n}\cdot\frac{\log\kappa_n}{\log n}\leq \frac{1}{\gamma}.
\end{equation}
Since $X_{\tau_{\kappa_n}}\cdot\vec{\ell}\leq X_n\cdot\vec{\ell}\leq X_{\tau_{\kappa_n+1}}\cdot\vec{\ell}$, we have that
\begin{align*}
|\Sc_{\kappa_n}-(X_n\cdot\vec{\ell}-vn)|
& \leq |X_{\tau_{\kappa_n+1}}\cdot\vec{\ell}-X_n\cdot\vec{\ell}|+v|\tau_{\kappa_n+1}-n|+|X_{\tau_{1}}\cdot\vec{\ell}-v\tau_{1}| \\
& \leq |X_{\tau_{\kappa_n+1}}\cdot\vec{\ell}-X_{\tau_{\kappa_n}}\cdot\vec{\ell}|+v|\tau_{\kappa_n+1}-\tau_{\kappa_n}|+|X_{\tau_{1}}\cdot\vec{\ell}-v\tau_{1}|.
\end{align*}
Let $\tilde{\varepsilon}>0$ be such that $(1/\gamma+\varepsilon)(1-\tilde{\varepsilon})\gamma>1$. Then, by Theorem \ref{t:up},
\[\sum_{k=1}^\infty\Pb_p(\tau_{k+1}-\tau_k>k^{1/\gamma+\varepsilon})\leq C\sum_{k=1}^\infty k^{-(1/\gamma+\varepsilon)(1-\tilde{\varepsilon})\gamma}<\infty.\]
Consequently, by the Borel-Cantelli lemma, we have that there are, $\Pb_p$-a.s., only finitely many $k$ such that $\tau_{k+1}-\tau_k>k^{1/\gamma+\varepsilon}$, and therefore only finitely many $n$ such that $\max_{k\leq n}(\tau_{k+1}-\tau_k)>n^{1/\gamma+\varepsilon}$. In particular, together with the convergence of $\kappa_n/n$, this implies that $v|\tau_{\kappa_n+1}-\tau_{\kappa_n}|$ is $\mathbb{P}_p$-a.s.\ bounded above by $cn^{1/\gamma+\varepsilon}$ for large $n$ (for some deterministic $c$). Similarly, using Theorem \ref{t:Xexp}, the same holds for $X_{\tau_{\kappa_n+1}}\cdot\vec{\ell}-X_{\tau_{\kappa_n}}\cdot\vec{\ell}$. In conjunction with \eqref{e:lskappa} and that $|X_{\tau_{1}}\cdot\vec{\ell}-v\tau_{1}|$ does not depend on $n$, we therefore find that
\[\limsup_{n\to\infty}\frac{\log|X_n\cdot\vec{\ell}-nv|}{\log n}\leq \frac{1}{\gamma},\]
as desired.

For the remaining claim, note that since $\Sc_n=(X_{\tau_{n+1}}-X_{\tau_1})\cdot\vec{\ell}-v(\tau_{n+1}-\tau_1)$ and $\Delta_{X_{\tau_{n+1}}\cdot\vec{\ell}}=\tau_{n+1}$,
\[\left|\Delta_{X_{\tau_{n+1}}\cdot\vec{\ell}}-v^{-1}X_{\tau_{n+1}}\cdot\vec{\ell}\right|=\left|\tau_{n+1}-v^{-1}X_{\tau_{n+1}}\cdot\vec{\ell}\right|=\left|v^{-1}\Sc_n + v^{-1}X_{\tau_1}\cdot\vec{\ell}-\tau_1\right|.\]
Applying \eqref{e:lsZ} and the fact that $X_{\tau_{n}}\cdot\vec{\ell}/n$ converges $\mathbb{P}_p$-a.s.\ to a limit in $(0,\infty)$, it follows that
\[\limsup_{s\rightarrow\infty}\frac{\log \left| \Delta_s-sv^{-1}\right|}{\log s}\geq\limsup_{n\rightarrow\infty}\left(\frac{\log \left|\Sc_n \right|}{\log n}\times \frac{\log n}{\log X_{\tau_{n+1}}\cdot\vec{\ell}}\right)=\frac{1}{\gamma}.\]
Moreover, for $X_{\tau_{n}}\cdot\vec{\ell}<s\leq X_{\tau_{n+1}}\cdot\vec{\ell}$, we have
\[0\leq \Delta_{X_{\tau_{n+1}}\cdot\vec{\ell}}-\Delta_s\leq \tau_{\vartheta_s+1}-\tau_{\vartheta_s}\leq \max_{k\leq n}\tau_{k+1}-\tau_k,\]
which is at most $n^{1/\gamma+\varepsilon}$ for all but finitely many $n$. In conjunction with the fact that there are, almost surely, only finitely many $n$ such that $X_{\tau_{n+1}}\cdot\vec{\ell}-X_{\tau_{n}}\cdot\vec{\ell}>n^{1/\gamma+\varepsilon}$, this ensures that
\[\limsup_{s\rightarrow\infty}\frac{\log \left| \Delta_s-sv^{-1}\right|}{\log s}\leq\frac{1}{\gamma}.\]
\end{proof}

\begin{proof}[Proof of \eqref{liminf}]
Again let $\Sc_n=\sum_{m=1}^nZ_m$. It then follows from \eqref{e:Zuplo} and Proposition \ref{p:limsup} that, $\Pb_p$-a.s.,
\[\liminf_{n\to\infty}\frac{\log_+|\Sc_n|}{\log n}=0.\]
Note again that $\Sc_n=(X_{\tau_{n+1}}-X_{\tau_1})\cdot\vec{\ell}-v(\tau_{n+1}-\tau_1)$, where $X_{\tau_1}\cdot\vec{\ell}$ and $\tau_1$ are $\Pb_p$-a.s.\ finite and do not depend on $n$. By \eqref{e:Ltn}, we then have that
\begin{align*}
\liminf_{n\to\infty}\frac{\log_+ |X_{\tau_{n+1}}\cdot\vec{\ell}-v\tau_{n+1}|}{\log \tau_{n+1}}=0, \qquad \Pb_p\text{-a.s.},
\end{align*}
and therefore
\begin{align*}
\liminf_{n\to\infty}\frac{\log_+ |X_{n}\cdot\vec{\ell}-vn|}{\log n}=0, \qquad \Pb_p\text{-a.s.}
\end{align*}
For the remaining claim, once more observing $\Delta_{X_{\tau_{n+1}}\cdot\vec{\ell}}=\tau_{n+1}$, we have that
\[\left|\Delta_{X_{\tau_{n+1}}\cdot\vec{\ell}}-v^{-1}X_{\tau_{n+1}}\cdot\vec{\ell}\right|=\left|\tau_{n+1}-v^{-1}X_{\tau_{n+1}}\cdot\vec{\ell}\right|=\left|v^{-1}\Sc_n + v^{-1}X_{\tau_1}\cdot\vec{\ell}-\tau_1\right|,\]
and so we obtain the result by applying the previous part of the proof again.
\end{proof}

\begin{proof}[Proof of \eqref{xlim}]
Let $\varepsilon>0$. By \eqref{limsup} we have that $\lim_{n\to\infty}\Pb(|X_n-nv|>n^{1/\gamma+\varepsilon})=0$, and therefore it remains to show that $\lim_{n\to\infty}\Pb(|X_n-nv|<n^{1/\gamma-\varepsilon})=0$.
Define $\overline{X}_n=\sup_{j\leq n}X_j\cdot \vec{\ell}$ to be the maximum distance reached by the walk in direction $\vec{\ell}$ up to time $n$. By Theorem \ref{t:Xexp}, we have
\begin{align*}
\Pb_p((X_{\tau_2}-X_{\tau_1})\cdot \vec{\ell} >A\log n)
& \leq Ce^{-cA\log n}
\end{align*}
for constants $c,C$. Choosing $A$ large and using Borel-Cantelli we thus have that
\[\Pb_p((X_{\tau_{k+1}}-X_{\tau_k})\cdot \vec{\ell}>A\log k\text{ i.o.})=0.\]
Since
\[0\leq\overline{X}_n-X_n\cdot\vec{\ell}\leq\max_{k\leq n}(X_{\tau_{k+1}}-X_{\tau_k})\cdot \vec{\ell},\]
it follows that
\[\Pb_p(\overline{X}_n-X_n\cdot\vec{\ell}>A\log n\text{ i.o.})=0.\]
To prove \eqref{xlim}, it therefore suffices to show that
\begin{align*}
\lim_{n\to\infty}\Pb(|\overline{X}_n-nv|<n^{1/\gamma-\varepsilon})=0.
\end{align*}
By the definition of $\overline{X}_n$ and $\Delta_s$ we have that $\{\Delta_s\leq n\}=\{\overline{X}_n\geq s\}$. In particular,
\[\{\overline{X}_n-nv< n^{1/\gamma-\varepsilon}\}  = \{\Delta_{nv+n^{1/\gamma-\varepsilon}}-(nv+n^{1/\gamma-\varepsilon})v^{-1}> -n^{1/\gamma-\varepsilon}v^{-1}\},\]
\[\{\overline{X}_n-nv\geq -n^{1/\gamma-\varepsilon}\}  = \{\Delta_{nv-n^{1/\gamma-\varepsilon}}-(nv-n^{1/\gamma-\varepsilon})v^{-1}\leq n^{1/\gamma-\varepsilon}v^{-1}\}.\]
Consequently, $\Pb(|\overline{X}_n-nv|<n^{1/\gamma-\varepsilon}) $ is at most
\begin{align*}
&\Pb_p\left(|\Delta_{nv-n^{1/\gamma-\varepsilon}}-(nv-n^{1/\gamma-\varepsilon})v^{-1}|\leq n^{1/\gamma-\varepsilon/2}v^{-1}\right)\\
& +
\Pb_p\left(\Delta_{nv+n^{1/\gamma-\varepsilon}}-\Delta_{nv-n^{1/\gamma-\varepsilon}}-2n^{1/\gamma-\varepsilon}v^{-1}\geq n^{1/\gamma-\varepsilon/2}v^{-1}(1-n^{-\varepsilon/2})\right).
\end{align*}
By \eqref{deltalim}, we have that the first term converges to $0$ as $n\to\infty$. Finally, since there can be at most $2n^{1/\gamma-\varepsilon}$ regenerations between levels $nv-n^{1/\gamma-\varepsilon}$ and $nv+n^{1/\gamma-\varepsilon}$, by applying Markov's inequality we find that
\begin{align*}
&\Pb_p\left(\Delta_{nv+n^{1/\gamma-\varepsilon}}-\Delta_{nv+n^{1/\gamma-\varepsilon}}-2n^{1/\gamma-\varepsilon}v^{-1}\geq n^{1/\gamma-\varepsilon/2}v^{-1}(1-n^{-\varepsilon/2})\right) \\
& \leq \Pb_p\left(\tau_{2n^{1/\gamma-\varepsilon}}\geq v^{-1}n^{\frac{1}{\gamma}-\varepsilon/2}(1+n^{-\varepsilon/2})\right) \\
& \leq \frac{2n^{1/\gamma-\varepsilon}\Eb_p[\tau_2-\tau_1]+\Eb_p[\tau_1]}{v^{-1}n^{\frac{1}{\gamma}-\varepsilon/2}(1+n^{-\varepsilon/2})} \\
& \leq Cn^{-\varepsilon/2},
\end{align*}
which converges to $0$ as $n\to\infty$.
\end{proof}

\section*{Acknowledgments}

This research was partially supported by JSPS KAKENHI, grant numbers 17H01093 and 19K03540 and by NUS grant R-146-000-260-114. It was initiated whilst AB was visiting DC at the Research Institute for Mathematical Sciences, an International Joint Usage/Research Center located in Kyoto University.

\bibliographystyle{amsplain}
\bibliography{brwsp}

\end{document}